# A GENERALIZATION OF THE KODAIRA VANISHING AND EMBEDDING THEOREM

YING ZHU

## Introduction

The Kodaira embedding theorem states that a compact complex manifold can be embedded in complex projective space $P^N$ for some $N$ if and only if the manifold admits a positive line bundle. More recently, Siu, Demailly and Riemanschneider have considered vanishing and embedding theorems for line bundles which are not strictly positive. More precisely, Riemenschneider [10] generalized the Kodaira embedding theorem to semi-positive line bundles for Kähler manifolds; Siu [13] generalized the Kodaira embedding theorem to semi-positive line bundles for non-Kähler manifolds, sloving the Grauert–Riemenschneider conjecture; Demailly [4] generalized the Kodaira embedding theorem to a line bundle whose curvature may be somewhat negative. In this paper, we adapt techniques developed by Elworthy–Rosenberg for vanishing theorems in Riemannian geometry to give a new extension of the Kodaira vanishing and embedding theorem.

Since the Kodaira vanishing theorem plays a key role in the proof of the Kodaira embedding theorem, we first prove the following generalization of the Kodaira vanishing theorem (cf. Theorem 25). Note that the hypotheses on the Riemannian geometry of the Kähler manifold are those which naturally occur in Cheeger–Gromov type compactness theorems.

**Theorem**: *Let $\mathcal{M} = \mathcal{M}(n, V, d, k_1, k_2)$ be the collection of pairs $(M, F)$ with $M$ a compact Kähler manifold of dimension $n$ with Ricci curvature greater than $k_1$, diameter less than $d$, and volume greater than $V$, and $F$ a Hermitian line bundle on $M$ with curvature $\Omega$ bounded by $k_2$. Let $(M, F) \in \mathcal{M}$ and $\lambda = \lambda(F \otimes K_M^*)$ be the lowest eigenvalue of the curvature $\Omega(F \otimes K_M^*)$. For $\alpha > 0$, there exists $a = a(\mathcal{M}, \alpha)$ such that if $\lambda > \alpha$ except on a set of volume less then $a$, then $H^q(M, \mathcal{O}(F)) = 0$, $\forall q > 0$.*

In another words, if $F$ is positive except on a set of small volume relative to the Riemannian geometry of $M$, then we still get a vanishing theorem. In particular, if $\lambda > 0$ everywhere, then $H^q(M, \mathcal{O}(F)) = 0$, $\forall q > 0$; this is the classical Kodaira vanishing theorem.

Via a general procedure for proving vanishing of cohomology given a Weitzenböck formula, we obtain the following vanishing theorem (cf. Theorem 15):

**Theorem**: *Let $M$ be a compact Kähler manifold and $F$ a Hermitian holomorphic line bundle on $M$. Let $\tau_1 \leq \tau_2 \leq \cdots \leq \tau_n$ be the eigenvalues of the curvature form of $F$. If there exists an $\epsilon > 0$ satisfying*

$$\| \min(\tau_1 + \tau_2 + \cdots + \tau_q - \tau_{q+1} - \cdots - \tau_n + q\omega_1 - \epsilon, 0) \|_{\frac{n}{2}} < \min(A^{-1}, \epsilon B^{-1})$$

*then $H^q(M, \mathcal{O}(F)) = 0$. Here $\omega_1$ is the smallest eigenvalue of the Ricci tensor and $A, B$ are positive constants which depend only on the dimension of $M$, an upper bound for the* *diameter of $M$, and lower bounds for the volume of $M$ and the Ricci curvature on $M$.*
In constrast to the previous theorem, this result does not imply the Kodaira vanishing theorem, so it is a new type of vanishing theorem.

The first vanishing theorem leads to our generalization (cf. Theorem 42) of the Kodaira embedding theorem for Kähler manifolds. The proof involves a careful analysis of the Riemannian geometry of a blow–up manifold.

**Theorem**: *Let $\mathcal{M} = \mathcal{M}(n, D_0, V_0, k_1, k_2, c_1, c_2)$ be the collection of pairs $(M, F)$ where*

(1) *$M$ is a compact Kähler manifold of complex dimension $n$ with $\operatorname{Vol}(M) \geq V_0$, $\operatorname{Diam}(M) \leq D_0$, $|\mathcal{R}ic(M)| \leq k_1$;*

(2) *On $M$ there exists an open set $U$ with volume at most $\frac{1}{2}V_0$ and a local coordinate $(U, \varphi)$ which satisfies (i) $\varphi(U)$ is a ball of radius 1 in $C^n$; (ii) $\|g(\frac{\partial}{\partial z^i}, \frac{\partial}{\partial \overline{z^j}})\|_{C^2} \leq c_1$; (iii) If $v = \sum_i a_i \frac{\partial}{\partial z^i} \in T_x^{1,0}M$, then $g(v,v) \geq c_2 \sum_i |a_i|^2$.*

(3) *$F$ is a holomorphic Hermitian line bundle on $M$ with curvature $|\Omega(F)| \leq k_2$.*

*Given $\alpha > 0$, there exists $c = c(\mathcal{M}, \alpha) > 0$ such that if $(M, F) \in \mathcal{M}$ and $\Omega(F) > \alpha$ except on a set in $M - U$ of volume less than $c$, then $M$ is projective algebraic.*

Again, the theorem states that if $M$ admits a line bundle which is positive except on a small set, them $M$ is projective algebraic. It is not hard to see that if $F$ is positive, then the hypotheses of the theorem are satisfied, which gives the Kodaira embedding theorem for Kähler manifolds. The proof involves showing that the Kähler manifold $M$ is Moishezon. Demailly [4] has also given an integral condition on the negative part of the curvature of a line bundle for a manifold to be Moishezon. However, his techniques are complex geometric, while our approch is based on a mixture of analysis and Riemannian geometry.

I would like to thank S. Rosenberg and S.-T. Yau for their help and encouragement.

## 1. Definitions, Notation and Basic Facts

Throughout this section, let $M$ be a compact Kähler manifold of complex dimension $n$ and $X$ a compact Riemannian manifold of dimension $m$. Let $B$ denote either $M$ or $X$ and $E$ a Riemannian (or Hermitian) vector bundle of rank $r$ on $B$ and $D$ a compatible connection on $E$.

(1) **Frame Field and Coframe Field**

$\{e_1, \ldots, e_m\}$ is called a frame field of $X$ near $P \in X$ if $e_1, \ldots, e_m$ are locally defined orthonormal vector fields on a neighborhood of $P$. $\{e_1, \ldots, e_m\}$ is called a normal frame field of $X$ near $P$ if $\{e_1, \ldots, e_m\}$ is a frame field of $X$ near $P$ and $\nabla_{e_i} e_j(P) = 0$ for all $i, j$. The dual frame field $\{e^1, \ldots, e^m\}$ of a (normal) frame field $\{e_1, \ldots, e_m\}$ of $X$ near $P$ is called a (normal) coframe field of $X$ near $P$. Similarly, $\{V_1, \ldots, V_n\}$ is called a frame field of type (1,0) on $M$ near $P$ if $V_1, \ldots, V_n$ are vector fields near $P$ of type (1,0) and $\langle V_i, \overline{V_j} \rangle = \delta_{ij}$ for $i, j$. The frame field $\{V_1, \ldots, V_n\}$ is normal if $(\nabla_{V_i} V^j)(P) = 0$, for all $i, j$. The dual frame $\{\omega^1, \ldots, \omega^n\}$ of a (normal) frame field $\{V_1, \ldots, V_n\}$ of type (1,0) on $M$ near $P$ is called a (normal) coframe field of type (1,0) on $M$ near $P$. $\{s_1, \ldots, s_r\}$ is called a frame of $E$ near $P$ if $s_1, \ldots, s_r$ are orthonormal local sections of $E$ near $P$.

(2) **Hermitian Metric, Inner Product**

A Hermitian metric on $M$ is a Riemannian metric $g$ invariant under the canonical complex structure $J$ on $M$. Let $\{V_1, \ldots, V_n\}$ be a frame field of type (1,0) on $M$ with dual frame $\{w^1, \ldots, w^n\}$. Write $V_k = \frac{1}{\sqrt{2}}(X_k - iY_k)$ and $\omega^k = \frac{1}{\sqrt{2}}(\alpha^k + i\beta^k)$ where $Y_k = JX_k$ and $\beta^k = J\alpha^k$ for $k = 1, \ldots, n$. Then $\{X_1, \ldots, X_n, Y_1, \ldots, Y_n\}$ is a frame field on $M$ (as a Riemannian manifold), and $\{\alpha^1, \ldots, \alpha^n, \beta^1, \ldots, \beta^n\}$ is the dual frame of $\{X_1, \ldots, X_n, Y_1, \ldots, Y_n\}$.
Any metric on $B$ induces an inner product $\langle \cdot, \cdot \rangle_B$ on $\Lambda^p T_x^* B$ for all $x \in B$. Given a metric $\langle \cdot, \cdot \rangle_E$, on the bundle $E$, we can define an inner product on $\Lambda^p T_x^* B \otimes E_x$ by

$$\langle \alpha \otimes s, \beta \otimes t \rangle = \langle \alpha, \beta \rangle_B \cdot \langle s, t \rangle_E,$$

where $\alpha, \beta \in \Lambda^p T_x^* B$ and $s, t \in E_x$. An inner product on $\Gamma(\Lambda^p T^* B \otimes E)$ is defined by

$$(\xi, \eta) = \int_B \langle \xi, \eta \rangle \, \mathrm{dvol}, \text{ for } \xi, \eta \in \Gamma(\Lambda^p T^* B \otimes E).$$

The pointwise norm of $\xi$ is defined by $|\xi|^2 = \langle \xi, \overline{\xi} \rangle$; the $L^2$–norm of $\xi$ is defined by $\|\xi\|^2 = (\xi, \overline{\xi})$. The Hodge $*$-operator is defined by $\xi \wedge *\eta = \langle \xi, \overline{\eta} \rangle \, \mathrm{dvol}$.

(3) $d, d^*, D, D^*$

Associated to a connection $D$ on $E$ we define $d^D : \Gamma(\Lambda^p T^* B \otimes E) \longrightarrow \Gamma(\Lambda^{p+1} T^* B \otimes E)$ by

$$d^D(\alpha \otimes s) = d\alpha \otimes s + (-1)^p \alpha \wedge Ds,$$

where $\alpha \in \Gamma(\Lambda^p T^* B)$ and $s \in \Gamma(E)$. We denote $d^D$ by $d$. Let $d^*$ be the adjoint of $d$ and $D^*$ the adjoint of $D$ with respect to the inner product on $\Gamma(\Lambda^p T^* B \otimes E)$. We define the **Riemannian Laplacian** on $\Gamma(\Lambda^p T^* B \otimes E)$ by $\Delta = dd^* + d^*d$. On a Hermitian holomorphic vector bundle, $d$ can be split into $d' + d''$ and $D$ can be split into $D' + D''$. Let $\Delta' = d'(d')^* + (d')^*d'$ and $\Delta'' = d''(d'')^* + (d'')^*d''$. Let $\overline{\partial} : \Gamma(\Lambda^{p,q} T^* X \otimes E) \longrightarrow \Gamma(\Lambda^{p,q+1} T^* X \otimes E)$ be defined by $\overline{\partial}(\alpha \otimes s) = \overline{\partial}\alpha \otimes s$, where $s$ is a holomorphic section. Let $\overline{\partial}^*$ be the adjoint of $\overline{\partial}$. We define the **complex Laplacian** $\square$ on $\Gamma(\Lambda^{p,q} T^* M \otimes E)$ by $\square = \overline{\partial}\overline{\partial}^* + \overline{\partial}^*\overline{\partial}$.

From the definition, we have $D = \sum_k e^k \otimes D_{e_k}, \quad d = \sum_k e^k \wedge D_{e_k}$, where $\{e_1, \ldots, e_m\}$ is a frame field on $X$ with dual frame $\{e^1, \ldots, e^m\}$. Let $\{V_1, \ldots, V_n\}, \{\omega^1, \ldots, \omega^n\}$ be as before. Then

$$D = \sum_{k=1}^n w^k \otimes D_{V_k} + \overline{w^k} \otimes D_{\overline{V_k}}, \quad d = \sum_{k=1}^n w^k \wedge D_{V_k} + \overline{w^k} \wedge D_{\overline{V_k}}.$$

**Remark:** We will frequently use this procedure to obtain formulas for the complex case from formulas for the real case. For example, an elementary computation shows that

$$d^* = -\sum_k i(e_k) D_{e_k},$$

where $i$ is the interior product. In the complex case, we have

(1) $$d^* = -\sum_k i(V_k) D_{\overline{V_k}} - \sum_k i(\overline{V_k}) D_{V_k}.$$

It is standard [3, §12, Chapter IV] that on manifolds

$$D^* D = -\sum_i D_{e_i} D_{e_i} + D_{D_{e_i} e_i},$$

and the proof carries over to the bundle case. In the complex case, we get

(2) $$D^* D = -\sum_k D_{V_k} D_{\overline{V_k}} - \sum_k D_{\overline{V_k}} D_{V_k}.$$

(4) **Eigenvalues of the Curvature Form of a Line Bundle**

Let $\Omega$ be the curvature form of $E$, defined by $\Omega = \overline{\partial}\partial \log |s|^2$, where $s$ is a holomorphic section of $E$. There exists a coframe field of type (1,0) $\{\xi^1, \cdots, \xi^n\}$ and real numbers $\lambda_1 \leq \lambda_2 \leq \cdots \leq \lambda_n$ such that
$$\Omega = \sum_j \lambda_j \xi^j \wedge \overline{\xi^j}.$$

We call $\{\lambda_i\}$ the eigenvalues of the curvature form $\Omega$. We say that $\Omega$ is bounded by $k$ if $|\lambda_i| \leq k$ for all $i$. This is equivalent to $|\Omega(v,v)| \leq k g(v,v)$ for all $v \in T^{1,0}M$. We also write $\Omega > \alpha$ if $\lambda_i > \alpha$ for all $i$.

## 2. A Weitzenböck Formula for the Complex Laplacian on a Hermitian holomorphic Vector Bundle

In this section, we are going to combine the Weitzenböck formula for the real Laplacian on a Riemannian vector bundle with the Kodaira–Nakano identity to obtain a new Weitzenböck formula for the complex Laplacian on a Hermitian holomorphic vector bundle. Let us first recall the Weitzenböck formula for the real Laplacian on a Riemannian bundle (cf. [7, Appendix]).

**Proposition 3.** *Let $X$ be a Riemannian manifold of dimension $m$, $E$ a Riemannian vector bundle over $X$ and $D$ a compatible connection on $E$. Then*
$$\Delta = D^*D - \sum_{j,k} e^j \wedge i(e_k) R_{e_j e_k},$$
*where $\{e_1, \ldots, e_m\}$ is a frame field on $X$ with its dual frame $\{e^1, \ldots, e^m\}$ and $R_{YZ} = D_Y D_Z - D_Z D_Y - D_{[Y,Z]}$ is the curvature of the connection $D$.*

For the rest of this section, let $M$ be a compact Kähler manifold of complex dimension $n$, $E$ a Hermitian holomorphic vector bundle over $M$ and $D$ a compatible connection on $E$. We can view $M$ as a real Riemannain manifold of dimension $2n$ and $E$ as a Riemannian vector bundle over $M$. In the notation of §1, we have
$$\begin{aligned}\Delta &= D^*D + \sum_{j,k} \alpha^j \wedge i(X_k) R_{X_k X_j} + \sum_{j,k} \beta^j \wedge i(Y_k) R_{Y_k Y_j} \\ &\quad + \sum_{j,k} \alpha^j \wedge i(Y_k) R_{Y_k X_j} + \sum_{j,k} \beta^j \wedge i(X_k) R_{X_k Y_j} \\ &= D^*D + \sum_{j,k} \omega^j \wedge i(V_k) R_{\overline{V_k} V_j} + \sum_{j,k} \overline{\omega^j} \wedge i(\overline{V_k}) R_{V_k \overline{V_j}} \\ &\quad - \sum_{j,k} \omega^j \wedge i(\overline{V_k}) R_{V_k V_j} - \sum_{j,k} \overline{\omega^j} \wedge i(V_k) R_{\overline{V_k V_j}}.\end{aligned}$$

Since $R_{JXJY} = R_{XY}$, the last two terms vanish. Set
$$(4) \qquad P = \sum_{j,k} \omega^j \wedge i(V_k) R_{\overline{V_k} V_j} + \sum_{j,k} \overline{\omega^j} \wedge i(\overline{V_k}) R_{V_k \overline{V_j}}.$$

Hence $\Delta = D^*D + P$. We define a linear mapping $L : \Gamma(\Lambda^{p,q} T^*M \otimes E) \longrightarrow \Gamma(\Lambda^{p+1,q+1} T^*M \otimes E)$ by $L(\alpha \otimes s) = (\psi \wedge \alpha) \otimes s$, where $\psi$ is the Kähler form on $M$. Let $L^*$ be the adjoint of $L$. The **Kodaira-Nakano Identity** [12, p. 16-17] is
$$(5) \qquad \Delta'' = \Delta' + i[d'd'' + d''d', L^*]$$

For $D$ be the cannonical connection, i.e. $d'' = \overline{\partial}$, we have $\square = \Delta''$, $\Omega = d^2 = d'd'' + d''d'$ and $\Delta = \Delta' + \Delta''$. Thus (5) becomes

$$(6) \qquad 2\square = \Delta + i[\Omega, L^*].$$

Combining the Weitzenböck formula $\Delta = D^*D + P$ and (6), we have

**Theorem 7 (Complex Weitzenböck formula).** *On a compact Kähler manifold, we have*

$$2\square = D^*D + R,$$

*where*

$$R = i[\Omega, L^*] + \sum_{j,k} \omega^j \wedge i(V_k) R_{\overline{V_k} V_j} + \sum_{j,k} \overline{\omega^j} \wedge i(\overline{V_k}) R_{V_k \overline{V_j}}.$$

**Remark**: $R$, as a linear endomorphism of $\Gamma(\Lambda^{p,q}T^*X \otimes E)$, involves two types of curvature, one from the manifold $X$ itself, and the other from the connection on the bundle $E$. We will analyze this more carefully in §3. From the computation in §3 we will see that $R$ is a symmetric endomorphism when $p = 0$ and $E$ is a line bundle.

## 3. Vanishing Theorem I for $H^q(M, \mathcal{O}(F))$

Let us first recall a general procedure for proving the vanishing of cohomology. Let $X$ be a compact Riemannian manifold and $E$ a Riemannian vector bundle over $X$. Let $D$ be a compatible connection on $E$. Assume that $\tilde{\Delta}$ is a naturally defined Laplacian having a Weitzenböck formula

$$\tilde{\Delta} = D^*D + R$$

where $R$ is a linear symmetric endomorphism of $E$.

Let $a_- = \min\{0, a\}$. Let $\rho(x)$ denote the lowest eigenvalue of $R_x$ on $E_x$, the fibre of $E$ at $x$, i.e.

$$\rho(x) = \inf\{\langle R_x(s), s\rangle_{E_x} : s \in E_x \quad \text{and} \quad \langle s, s\rangle_{E_x} = 1\}.$$

The main result we are going to use in this section is the following theorem of Rosenberg–Yang [11]:

**Proposition 8.** *If there exists $a > 0$ such that*

$$\|(\rho - a)_-\|_{\frac{n}{2}} < min(A^{-1}, aB^{-1}),$$

*then the operator $\tilde{\Delta}$ is strictly positive on $E$. Here $A, B$ are positive constants which depend only on the dimension of $X$, an upper bound for the diameter of $X$, and lower bounds for the volume of $X$ and the Ricci curvature on $X$.*

For the rest of this section, let $M$ be a compact Kähler manifold and $F$ a Hermitian holomorphic line bundle over $M$. Let $E = \Lambda^{(0,q)}T^*M \otimes F$ and $D$ the canonical connection on $F$. We extend $D$ to be a connection on $E$ by tensoring it with the Levi-Civita connection $\nabla$ on $M$. Let $\tilde{\Delta} = 2\square$. We have from Theorem 7 that $\tilde{\Delta} = 2\square = D^*D + i[\Omega, L^*] + P$, where $P$ is as in (4).

Since $H^q(M, \mathcal{O}(F)) = \text{Ker}\,\square$, the positivity of $\square$ implies the vanishing of $H^q(M, \mathcal{O}(F))$. In order to get a sufficient condition for the vanishing of $H^q(M, \mathcal{O}(F))$ by Proposition 8, we need to compute $\rho$ for $i[\Omega, L^*] + P$. We separate the computation into two parts:

**Part 1**: Computation of $\langle i[\Omega, L^*]\Phi, \overline{\Phi}\rangle$ where $\Phi \in \Gamma(E) = \Gamma(\Lambda^{0,q}T^*M \otimes F)$. Write $\Phi = e \cdot \phi$, with $e$ a frame of $F$, and $\phi$ a (0,q)-form on $M$. There exists a coframe field $\{\eta^1, \cdots, \eta^n\}$ of type (1,0) on $M$ such that $\Omega = \sum_k \tau_k(x)\eta^k \wedge \overline{\eta^k}$, where $\tau_1 \leq \tau_2 \leq \cdots \leq \tau_n$ are the eigenvalues of $\Omega$. The Kähler form on $M$ is $\psi = i\sum_k \eta^k \wedge \overline{\eta^k}$. Write $\phi = \sum_I a_I \overline{\eta^I}$, where $\eta^I = \eta^{i_1} \wedge \cdots \wedge \eta^{i_q}$ and $I = (i_1, \cdots, i_q)$ ranges over all multi-indices of length $q$ with $i_1 < i_2 < \cdots < i_q$. Then we have

$$\begin{aligned}
\langle i[\Omega, L^*]\Phi , \overline{\Phi}\rangle &= -i\langle L^*(\Omega \wedge \Phi) , \overline{\Phi}\rangle \\
&= -i\langle \Omega \wedge \Phi , \overline{\psi} \wedge \overline{\Phi}\rangle \\
&= \langle \sum_I a_I \sum_{k \notin I} \tau_k(x)\eta^k \wedge \overline{\eta^k} \wedge \overline{\eta^I} , \sum_{l,J} \overline{a_J}\eta^l \wedge \overline{\eta^l} \wedge \eta^J\rangle \\
&= -\sum_I (\sum_{k \notin I} \tau_k(x))|a_I|^2.
\end{aligned}$$

**Part 2**: Computation of $\langle P\Phi, \overline{\Phi}\rangle$. Since the first term in $P$ vanishes on $\Gamma(\Lambda^{0,q}T^*M \otimes F)$, we have

$$(9) \qquad P\Phi = \sum_{j,k} \overline{\xi^j} \wedge i(\overline{V_k})R_{V_k\overline{V_j}}\Phi,$$

where $\{V_1, \ldots, V_n\}$ is a frame field of type (1,0) on $M$ with dual frame $\{\xi^1, \ldots, \xi^n\}$. Note that $R_{XY}\Phi = 2\Omega(X,Y)\Phi + e \cdot R'_{XY}\phi, \forall X,Y \in T_pM$, where $e$ is a frame of $F$, $\Phi = e \cdot \phi \in \Gamma(\Lambda^{p,q}T^*M \otimes F)$, and $R'_{XY}$ is the curvature on the Riemannian manifold $M$, and $\Omega$ is the curvature of the line bundle $F$. (This follows from the identity $\Omega_{\Lambda^{p,q}T^*M \otimes F} = \Omega_{\Lambda^{p,q}T^*M} \otimes I + I \otimes \Omega_F$.) Thus (9) becomes

$$(10) \qquad P\Phi = 2\sum_{j,k} \overline{\xi^j} \wedge i(\overline{V_k})\Omega(V_k,\overline{V_j})\Phi + e \cdot \sum_{j,k} \overline{\xi^j} \wedge i(\overline{V_k})R'_{V_k\overline{V_j}}\phi.$$

**Lemma 11.**

$$e \cdot \sum_{j,k} \overline{\xi^j} \wedge i(\overline{V_k})R'_{V_k\overline{V_j}}\phi = 2\sum_{j,k} \overline{\xi^j} \wedge i(\overline{V_k})\Omega_{K_M^*}(V_k,\overline{V_j})\Phi,$$

where $K_M$ is the canonical line bundle of $M$.

**Proof:** We may assume that $\{V_1, \ldots, V_n\}$ and $\{\xi^1, \ldots, \xi^n\}$ are chosen such that $\Omega_{K_M^*} = \sum_j \omega_j \xi^j \wedge \overline{\xi^j}$. The Ricci curvature tensor and the curvature of the canonical line bundle are related by $\mathcal{R}ic(V_i, \overline{V_j}) = \omega_i \delta_{ij}$. Since $\square = \overline{\square}$ on Kähler manifolds, we have

$$(12) \qquad \sum_{j,k} \overline{\xi^j} \wedge i(\overline{V_k})R'_{V_k\overline{V_j}} = \sum_k R'_{V_k\overline{V_k}} \qquad \text{on} \quad \Gamma(\Lambda^{0,q}T^*M).$$

From [16, p. 952-953] we know

$$(13) \qquad \sum_k R'_{V_k\overline{V_k}}\xi^j = -\mathcal{R}ic(V_k,\overline{V_k})\xi^j = -2\sum_{j,k} \overline{\xi^j} \wedge i(\overline{V_k})\Omega_{K_M^*}(V_k,\overline{V_j})\xi^j.$$

From (12, 13) we get the result.　　　　　　　　　　　　　　　　QED

By Lemma 11, (10) becomes

$$(14) \qquad P\Phi = 2\sum_{j,k} \overline{\xi^j} \wedge i(\overline{V_k})\Omega(V_k,\overline{V_j})\Phi + 2\sum_{j,k} \overline{\xi^j} \wedge i(\overline{V_k})\Omega_{K_M^*}(V_k,\overline{V_j})\Phi.$$

For $\Omega_{K_M^*} = \sum_k \omega_k \xi^k \wedge \overline{\xi^k}, \Omega = \sum_k \tau_k \eta^k \wedge \overline{\eta^k}, \phi = \sum_I b_I \overline{\xi}^I$, we have

$$\begin{aligned}\langle P\Phi, \overline{\Phi}\rangle &= 2\langle \sum_{j,k} \overline{\eta^j} \wedge i(\overline{W_k})\Omega(W_k, \overline{W_j})\Phi, \overline{\Phi}\rangle + 2\langle \sum_{j,k} \overline{\xi^j} \wedge i(\overline{V_k})\Omega_{K_M^*}(V_k, \overline{V_j})\Phi, \overline{\Phi}\rangle \\ &= \sum_I (\sum_{k \in I} \tau_k(x))|a_I|^2 + \sum_I (\sum_{k \in I} \omega_k(x))|b_I|^2,\end{aligned}$$

where $\{W_1, \ldots, W_n\}$ is the dual frame of $\{\eta^1, \ldots, \eta^n\}$. Thus

$$\begin{aligned}\rho &= \inf_\Phi \{\langle i[\Omega, L^*]\Phi, \overline{\Phi}\rangle + \langle P\Phi, \overline{\Phi}\rangle : \Phi \in \Gamma(\Lambda^{0,q}T^*M \otimes F), \langle \Phi, \overline{\Phi}\rangle = 1\} \\ &= \inf\{\sum_I (\sum_{k \in I} \tau_k - \sum_{k \notin I} \tau_k)|a_I|^2 + \sum_I (\sum_{k \in I} \omega_k)|b_I|^2 : \sum_I |a_I|^2 = 1, \sum_I |b_I|^2 = 1\} \\ &\geq \tau_1 + \tau_2 + \cdots + \tau_q - \tau_{q+1} - \cdots - \tau_n + q\omega_1.\end{aligned}$$

Applying Proposition 8, we get

**Theorem 15.** *Let $M$ be a compact Kähler manifold and $F$ a Hermitian holomorphic line bundle on $M$ with curvature form $\Omega$. If there exists $\epsilon > 0$ satisfying*

$$\|(\tau_1 + \tau_2 + \cdots + \tau_q - \tau_{q+1} - \cdots - \tau_n + q\omega_1 - \epsilon)_-\|_{\frac{n}{2}} < \min(A^{-1}, \epsilon B^{-1}),$$

*then $H^q(M, \mathcal{O}(F)) = 0$ for $q > 0$. Here $\omega_1$ is the smallest eigenvalue of the Ricci tensor, and $\tau_1 \leq \tau_2 \leq \cdots \leq \tau_n$ are the eigenvalues of $\Omega$, $A, B$ are positive constants which depend only on the dimension of $M$, an upper bound for the diameter of $M$, and lower bounds for the volume of $M$ and the Ricci curvature on $M$.*

Notice that since the first chern class $c_1(M) \geq 0$ if and only if $\omega_1 \geq 0$, we have

**Corollary 16.** *Let $M$ be a compact Kähler manifold with $c_1(M) > 0$ (as a (1,1) form). Let $F$ be a Hermitian holomorphic line bundle on $M$ with curvature form $\Omega$. If there exists $\epsilon > 0$ satisfying*

$$\|(\tau_1 + \tau_2 + \cdots + \tau_q - \tau_{q+1} - \cdots - \tau_n - \epsilon)_-\|_{\frac{n}{2}} < \min(A^{-1}, \epsilon B^{-1})$$

*then $H^q(M, \mathcal{O}(F)) = 0$, where $\tau_i, A, B$ are as in Theorem 15.*

**Corollary 17.** *Let $M$ be a compact Kähler manifold, $F$ a Hermitian holomorphic line bundle on $M$ with curvature $\Omega$. Let $\tau_1 \leq \tau_2 \cdots \leq \tau_n$ be the eigenvalues of $\Omega$. If $\tau_1 + \cdots + \tau_q > \tau_{q+1} + \cdots + \tau_n$, then for all $k$ sufficiently large, we have $H^q(M, \mathcal{O}(F^k)) = 0$.*

**Corollary 18.** *Let $M$ be a compact Kähler manifold. If the Ricci curvature is positive, then $H^{0,q}(M, \mathbb{C}) = 0$ for all $q > 0$.*

Corollary 18, which is a direct consequence of the Kodaira vanishing theorem, can be generalized as in the following

**Corollary 19.** *Let $M$ be a compact Kähler manifold. Let $\omega$ be the smallest eigenvalue of the Ricci curvature. If there exists $\epsilon > 0$ such that*

$$\|(\omega - \epsilon)_-\|_{\frac{n}{2}} \leq \min(A^{-1}, \epsilon B^{-1}),$$

*then $H^{0,q}(M, \mathbb{C}) = 0$ for all $q > 0$, where $A, B$ are as in Theorem 15.*

## 4. Vanishing Theorem II for $H^{0,q}(M,F)$

Throughout this section, let $M$ be a compact Kähler manifold and $\nabla$ the Levi-Civita connection on $M$. Let $F$ be a Hermitian holomorphic line bundle on $M$ and $D$ the canonical connection on $F$. The Kodaira vanishing theorem says that if $F \otimes K_M^*$ is positive, then $H^q(M, \mathcal{O}(F))$ vanishes for $q > 0$. Here the positivity of a bundle means the positivity of the lowest eigenvalue of the curvature of the bundle. We first rewrite Kodaria's original proof [6, §3.7] in a coordinate free manner, as we will generalize these techniques below.

**Theorem 20.** *Let $\lambda = \lambda(F \otimes K_M^*)$ be the lowest eigenvalue of the curvature $\Omega(F \otimes K_M^*)$. If $\Phi \in \Gamma(\Lambda^{0,q} T^*M \otimes F)$ has $\Box \Phi = 0$, then $\int_M \lambda |\Phi|^2 \, dvol \leq 0$, where $q > 0$.*

**Corollary 21 (Kodaira Vanishing Theorem).** *If $F \otimes K_M^*$ is positive, then $H^q(M, \mathcal{O}(F)) = 0$ for all $q > 0$.*

Note that this theorem immediately implies a generalization of the Kodaira vanishing theorem due to Shiffman-Sommese [12, p. 34-36], namely that $F \otimes K_M^*$ semipositive everywhere and positive somewhere implies $H^q(M, \mathcal{O}(F)) = 0$ for all $q > 0$.

**Proof of Theorem 20**: Assume $\Phi \in \Gamma(\Lambda^{0,q} T^*M \otimes F)$ and $\Box \Phi = 0$. Let $\beta = \langle D''\Phi, \overline{\Phi} \rangle$. Since $\beta$ is a (0,1) form on $M$, we have $\int_M d * \beta = 0$ by Stokes' theorem. Let $\{V_1, \ldots, V_n\}$ be a normal frame field of type (1,0) on $M$ with dual frame $\{\xi^1, \ldots, \xi^n\}$. By (1) we have

$$\begin{aligned}
d^*\beta = \overline{\partial}^*\beta &= -\sum_k i(\overline{V_k})(\nabla_{V_k}\beta) = -\sum_k (\nabla_{V_k}\beta)(\overline{V_k}) \\
&= -\sum_k V_k \cdot \beta(\overline{V_k}) = -\sum_k V_k \cdot \langle D_{\overline{V_k}}\Phi, \overline{\Phi} \rangle \\
&= -\sum_k \langle D_{V_k} D_{\overline{V_k}}\Phi, \overline{\Phi} \rangle - \sum_k \langle D_{\overline{V_k}}\Phi, \overline{D_{\overline{V_k}}\Phi} \rangle.
\end{aligned}$$

Since $\sum_k \langle D_{\overline{V_k}}\Phi, \overline{D_{\overline{V_k}}\Phi} \rangle \geq 0$, we have

$$(22) \qquad \int_M \langle \sum_k D_{V_k} D_{\overline{V_k}}\Phi, \overline{\Phi} \rangle \, dvol \leq 0.$$

Since $D$ is the canonical connection, $\overline{\partial} = d'' = \sum_j \overline{\xi^j} \wedge D_{\overline{V_j}}$ and $\overline{\partial}^* = -\sum_k i(\overline{V_k}) D_{V_k}$. Then the complex analogue of Proposition 3 gives $\Box = \sum_k D_{V_k} D_{\overline{V_k}} - \sum_{j,k} \overline{\xi^j} \wedge i(\overline{V_k}) R_{V_k \overline{V_j}}$. Thus from (22), we have $\int_M \langle \sum_{j,k} \overline{\xi^j} \wedge i(\overline{V_k}) R_{V_k \overline{V_j}} \Phi, \overline{\Phi} \rangle \, dvol \leq 0$. Since $\Phi$ is of type (0,q), by (14) we have

$$(23) \qquad \sum_{j,k} \overline{\xi^j} \wedge i(\overline{V_k}) R_{V_k \overline{V_j}} \Phi = 2 \sum_{j,k} \overline{\xi^j} \wedge i(\overline{V_k})(\Omega + \Omega_{K_M^*})(V_k, \overline{V_j}) \Phi.$$

Since both sides of the equation (23) are independent of the choice of frame field, we may assume that $\{V_1, \ldots, V_n\}$ and $\{\xi^1, \ldots, \xi^n\}$ are chosen such that $\Omega + \Omega_{K_M^*} = \sum_j \lambda_j \xi^j \wedge \overline{\xi^j}$, with $\lambda_1 \leq \lambda_2 \leq \cdots \leq \lambda_n$ the eigenvalues of $\Omega + \Omega_{K_M^*}$. Write $\Phi = e \cdot \phi$, $\phi = \sum_I a_I \overline{\xi}^I$ as usual. (23) becomes

$$\begin{aligned}
\sum_{j,k} \overline{\xi^j} \wedge i(\overline{V_k}) R_{V_k \overline{V_j}} \Phi &= e \cdot \sum_k \overline{\xi^k} \wedge i(\overline{V_k}) \lambda_k (\sum_I a_I \overline{\xi}^I) \\
&= e \cdot \sum_I (\sum_{k \in I} \lambda_k) a_I \overline{\xi}^I.
\end{aligned}$$

Thus

$$\begin{aligned}
0 &\geq \int_M \langle \sum_{j,k} \overline{\xi^j} \wedge i(\overline{V_k}) R_{V_k \overline{V_j}} \Phi \,, \overline{\Phi} \rangle \, \text{dvol} \\
&= \int_M \langle \sum_I (\sum_{k \in I} \lambda_k) a_I \overline{\xi}^I \,, \sum_I \overline{a_I} \xi^I \rangle \text{dvol} \\
&\geq \int_M \sum_I |a_I|^2 (\lambda_1 + \cdots + \lambda_q) \, \text{dvol} \\
&\geq q \int_M |\Phi|^2 \lambda_1 \, \text{dvol} \,.
\end{aligned}$$

Since $\lambda_1 = \lambda = \lambda\left(\Omega(F \otimes K_M^*)\right)$, we are done. QED

From the proof of Theorem 20, we easily get the following corollary:

**Corollary 24.** *(cf. [12, p. 34-36]) Let $\lambda_1 \leq \lambda_2 \leq \cdots \leq \lambda_n$ be the eigenvalues of the curvature $\Omega(F \otimes K_M^*)$. If $\lambda_1 + \cdots + \lambda_q \geq 0$ everywhere and is strictly positive somewhere, then $H^q(M, \mathcal{O}(F)) = 0$.*

Now we sharpen Theorem 20 by proving a vanishing theorem when $\lambda$ is negative on a set of small volume.

**Theorem 25.** *Let $\mathcal{M} = \mathcal{M}(n, V, d, k_1, k_2)$ be the collection of pairs $(M, F)$ with $M$ a compact Kähler manifold of dimension $n$ with Ricci curvature greater than $-k_1$ where $k_1 \geq 0$, diameter less than $d$, and volume greater than $V$ and $F$ a Hermitian holomorphic line bundle on $M$ with curvature $\Omega$ bounded by $k_2$. Let $\lambda = \lambda\left(F \otimes K_M^*\right)$ be the lowest eigenvalue of the curvature $\Omega\left(F \otimes K_M^*\right)$. For $\alpha > 0$, there exists $a = a(\mathcal{M}, \alpha)$ such that for any $(M, F) \in \mathcal{M}$ if $\lambda > \alpha$ except on a set of volume less then $a$, then $H^q(M, \mathcal{O}(F)) = 0, \forall q > 0$.*

It is easy to see that Theorem 25 generalizes the Kodaira vanishing theorem. The proof of Theorem 25 depends on Li's $L^\infty$ estimates for eigenfunctions of the Laplacian [9] and Elworthy–Rosenberg's estimates for regions of negative curvature [5]. We begin with a preparatory lemma:

**Lemma 26.** *Let $M$ be a compact Kähler manifold of dimension $n$ with Ricci curvature greater than $-k_1$, and $F$ a Hermitian holomorphic line bundle on $M$ with curvature bounded by $k_2$. If $\Phi \in \Gamma\left(\Lambda^{0,q} T^*M \otimes F\right)$ has $\Box \Phi = 0$, then*

$$\Delta |\Phi| \leq (qk_1 + nk_2) |\Phi|.$$

**Proof:** By (2) and the fact that $\Delta = D^*D$ on $C^\infty(M)$, it is not hard to show that

(27) $$\Delta \langle \Phi \,, \overline{\Phi} \rangle = -2 \langle D\Phi \,, \overline{D\Phi} \rangle + \langle D^*D\Phi \,, \overline{\Phi} \rangle + \langle \Phi \,, D^*D\overline{\Phi} \rangle.$$

From the Weitzenböck formula $2\Box = D^*D + R$ on $\Gamma\left(\Lambda^{p,q} T^*M \otimes F\right)$ and the assumption $\Box \Phi = 0$, we get $D^*D\Phi = -R\Phi$. From the computation in §3, we see that $\langle R\Phi \,, \overline{\Phi} \rangle$ is real for $\Phi \in \Gamma\left(\Lambda^{0,q} T^*M \otimes F\right)$. Hence $\langle D^*D\Phi \,, \overline{\Phi} \rangle$ is real and $\langle D^*D\Phi \,, \overline{\Phi} \rangle = \langle \Phi \,, D^*D\overline{\Phi} \rangle$. By (27) we have

(28) $$\Delta \langle \Phi \,, \overline{\Phi} \rangle = -2 \langle D\Phi \,, \overline{D\Phi} \rangle + 2 \langle D^*D\Phi \,, \overline{\Phi} \rangle.$$

On the other hand,

(29) $$\Delta \langle \Phi , \overline{\Phi} \rangle = \Delta |\Phi|^2 = 2|\Phi|\Delta|\Phi| - 2|\nabla|\Phi||^2.$$

By (28), (29) and the first Kato inequality (cf. [1] [2]): $|\nabla|\Phi||^2 \leq |D\Phi|^2$, we have

$$|\Phi|\Delta|\Phi| \leq \langle D^*D\Phi , \overline{\Phi} \rangle = \langle -R\Phi , \overline{\Phi} \rangle \leq -\rho|\Phi|^2,$$

where $\rho \geq \tau_1 + \tau_2 + \cdots + \tau_q - \tau_{q+1} - \cdots - \tau_n + q\omega_1 \geq -qk_1 - nk_2$ from the computations in §3. QED

**Proof of Theorem 25**: Assume $\Phi \in \Gamma(\Lambda^{0,q}T^*M \otimes F)$ and $\Box\Phi = 0$. Applying Lemma 26 and the proof of Theorem 7 in [9], we get

$$\|\Phi\|_\infty^2 \leq C_1 \|\Phi\|_2^2,$$

where $C_1 = C_1(\mathcal{M}) > 0$ is a constant depends only on $n, d, V, k_1, k_2$. Assume $\|\Phi\|_2^2 = 1$ and let $f$ denote $|\Phi|$. Thus we have $\|f\|_\infty^2 \leq C_1$, $\int_M f^2 \, d\text{vol} = 1$. From the proof of Proposition 1.2 in [5], we know that for $\alpha > 0$, there exists $a = a(C_1, \mathcal{M}, \alpha) = a(\mathcal{M}, \alpha)$ such that if $\lambda > \alpha$ except on a set of volume less then $a$, then

$$\int_M \lambda f^2 \, d\text{vol} > 0,$$

which contradicts Theorem 20. Thus $\Phi = 0$ and $H^q(M, \mathcal{O}(F))$ vanishes. QED

**Remark**: As in Corollary 24 we may replace $\lambda$ in Theorem 25 by $\lambda_1 + \cdots + \lambda_q$.

## 5. Moisezon Manifolds, Vanishing Theorems and the geometry of blow–up manifolds

A Moisezon space is a compact complex space of complex dimension $n$ and transcendence degree $n$. In particular, any algebraic variety is a Moisezon space. Wells [15] gives a simple criterion for a compact complex space to be Moisezon:

**Lemma 30.** *Let $f : X \longrightarrow \mathbb{P}^N$ be a meromorphic mapping, where $X$ is an irreducible complex space of dimension $\leq N$, $f$ is a holomorphic mapping on $X - S$ for $S$ a subvariety of $X$. If there is a point $x_o \in X - S$ such that*

*a) $x_o$ is not a singular point of $X$,*

*b) $f$ has maximal rank at $x_o$,*

*then $X$ is a Moisezon space.*

Let $m_p$ denote the sheaf of the holomorphic functions which vanish at $p$. Combining Lemma 30 and the proof of the Kodaira embedding theorem in ([14, p. 234-239], we get the following theorem:

**Theorem 31.** *Let $M$ be a compact complex manifold. If there is a holomorphic line bundle $F$ on $M$ and a point $p \in M$ such that $H^1(M, \mathcal{O}(F) \otimes m_p^\mu) = 0$ for $\mu = 1, 2$, then $M$ is a Moisezon manifold.*

From [15] we know that a Moisezon Kähler manifold is actually projective algebraic, so we have the following

**Corollary 32.** *Let $M$ be a compact Kähler manifold. If there exists a holomorphic line bundle $F$ on $M$ and a point $p \in M$ such that $H^1(M, \mathcal{O}(F) \otimes m_p^\mu) = 0$ for $\mu = 1, 2$, then $M$ can be embedded in a projective space.*

Since $m_p$ is not a locally free sheaf, we cannot apply our generalization of the Kodaira vanishing theorem. However we can apply it on the blow–up manifold, $\tilde{M}_p$, whose definition we now recall. Let $U$ be a neighborhood of the point $p$ with coordinates $z = (z_1, \cdots, z_n)$, where $z = 0$ corresponds to $p$. Let $(l_1, \cdots, l_n)$ be homogeneous coordinates for $\mathbb{P}^{n-1}$. Let $W = \{(z, l) = (z_1, \cdots, z_n, [l_1, \cdots, l_n]) \mid z_i l_j = z_j l_i \text{ for } i, j = 1, \cdots, n\}$. Let $\pi : W \to U, \pi(z, l) = z$ be the projection on the first factor. We define the blow–up $\tilde{M}_p$ of $M$ at $p$ to be the complex manifold $\tilde{M}_p = (M - p) \cup_\pi W$. For $S \stackrel{def}{=} \pi^{-1}(p) = \{0\} \times \mathbb{P}^{n-1}$, $\pi$ is a biholomorphic map from $\tilde{M}_p - S$ onto $M - p$. Let $L_p$ be the holomorphic line bundle on $\tilde{M}_p$ associated to the divisor $S$. From the proof of Kodaira embedding theorem in [14, p. 234-239] we get the following

**Lemma 33.** *For any $\mu \in \mathbb{Z}^+$, if $H^1(\tilde{M}_p, \pi^* F \otimes L_p^{*\mu}) = 0$, then $H^1(M, \mathcal{O}(F) \otimes m_p^\mu) = 0$.*

By Lemma 33, the vanishing of $H^1(\tilde{M}_p, \pi^* F \otimes L_p^{*\mu})$ for $\mu = 1, 2$ implies that $M$ is projective algebraic. In order to apply our generalization of the Kodaira vanishing theorem, we need to understand the geometry of $\tilde{M}_p$. We will devote the rest of this section to estimating diameter, volume and curvature bounds of the manifold $\tilde{M}_p$ in terms of those quantities on the base manifold $M$.

Let $g$ be a Kähler metric on $M$. Let $\tilde{g}$ be a metric on $\tilde{M}_p$ defined by

$$\tilde{g} = A\pi^* g + \overline{\partial}\partial(\log|z|^2 \cdot \pi^*\rho),$$

where $|z|^2 = |z_1|^2 + \cdots + |z_n|^2$, $\rho$ is a positive bump function which is 1 near $p$ and 0 outside $U$ and $A$ is some positive constant. Since $\pi^* g$ is positive definite where $\overline{\partial}\partial(\log|z|^2 \cdot \pi^*\rho)$ may be negative, we easily can choose $A$ big enough so that $\tilde{g}$ is positive definite. Since $d\tilde{g} = 0$, $\tilde{g}$ is a Kähler metric on $\tilde{M}_p$. In order to estimate the volume, diameter and Ricci curvature of $\tilde{M}_p$, we need to estimate $A$.

**Remark**: For conveniencem from now on we just write $g(v, w)$ for $g(v, \overline{w})$.

**Definition**: A local holomorphic coordinate $(U, \varphi) = (U, z_1, \ldots, z_n)$ on $M$ is said to be **of class $\mathcal{C}(c_1, c_2)$** if

(1)  $\varphi(U)$ is a ball of radius 1 in $\mathbb{C}^n$,
(2)  $\|g_{i\bar{j}}\|_{C^2} \leq c_1$,
(3)  $g \geq c_2 g_\varphi$.

Here $g_\varphi$ is the Euclidean metric on $U$ induced by $\varphi$, so (3) means that if $v = \sum_i a_i \frac{\partial}{\partial z_i} \in T_x^{1,0} M$, where $x \in U$, then $g(v, v) \geq c_2 \sum_i |a_i|^2$.

**Note**: For any coordinate $(U, \varphi)$ there always exists $c_1, c_2 > 0$ such that $(U, \varphi)$ is of class $\mathcal{C}(c_1, c_2)$.

In the construction of $\tilde{M}_p$, let $(U, \varphi)$ be a fixed chart centered at $p$ of class $\mathcal{C}(c_1, c_2)$. Let $\phi$ be a fixed cut–off function on $\mathbb{R}$ such that $\phi(x) = 0$, if $|x| \geq \frac{3}{4}$; and $\phi(x) = 1$, if $|x| \leq \frac{1}{2}$. Let $\rho(z) = \phi(|z|^2)$, and $U' = \varphi^{-1}(B(\frac{1}{2}))$, where $B(\frac{1}{2})$ is the ball of radius $\frac{1}{2}$ centered at the origin in $\mathbb{C}^n$. Then $\rho$ is a cut-off function on $M$ such that

$$\rho(x) = \begin{cases} 0 & x \notin U \\ 1 & x \in U'. \end{cases}$$

Denote $\frac{\partial}{\partial z_i}$ by $\partial_i$, and $\frac{\partial}{\partial \overline{z_j}}$ by $\partial_{\bar{j}}$. We have $\partial_i \rho = d\phi \cdot \overline{z_i}$, $\partial_{\bar{i}} \rho = d\phi \cdot z_i$, $\partial_i \partial_{\bar{j}} \rho = d^2\phi \cdot \overline{z_i} z_j$, etc., so there exists a constant $C(n)$ depending only on $n$ such that

$$\|\rho\|_{C^4} \leq C(n).$$

Let $\tilde{U} = \pi^{-1}(U)$ and $\tilde{U}' = \pi^{-1}(U')$. Given coordinates $(U, \varphi) = (U, z^1, \cdots, z^n)$ on $M$ centered at $p$, we can define two different coordinates which will be useful in later computations as follows:

**Definition: The pull–back coordinates** $(\tilde{U} - \tilde{U}', \varphi') = (\tilde{U} - \tilde{U}', z'_1, \ldots, z'_n)$ are defined by

$$(z_1, \ldots, z_n, [l_1, \ldots, l_n]) \xrightarrow{\varphi'} (z_1, \ldots, z_n) \text{ i.e. } z'_i = z_i, \text{ for all } i.$$

**Definition: The new coordinates** on $\tilde{U}'$. Let $\tilde{U}'_i = \{(z, l) \in \tilde{U}' \,|\, l_i \neq 0 \text{ and } |\frac{l_j}{l_i}| < 10, \text{ if } j \neq i\}$. Then $\{\tilde{U}'_i\}$ cover $\tilde{U}'$. Let $(\tilde{U}'_i, \tilde{\varphi}_i) = (\tilde{U}'_i, \tilde{z}_1, \ldots, \tilde{z}_n)$ be the coordinates defined by

$$\tilde{z}_j = \frac{l_j}{l_i} = \frac{z_j}{z_i} \text{ for } j \neq i; \quad \tilde{z}_i = z_i.$$

**Remark**: In the rest of this paper, we make the conventions that $C(x, y, z, \cdots)$ denotes a constant depending only on $x, y, z, \cdots$, and that the value of $C(x, y, z, \cdots)$ may change from line to line in a proof.

**Lemma 34.** *Let $\tilde{M}_p$ be the blow–up manifold of $M$ at $p$, $(U, \varphi)$ a holomorphic coordinate centered at $p$ of class $\mathcal{C}(c_1, c_2)$. There exists a positive constant $A(c_2, n)$ such that*
  (1) *$\tilde{g} = A(c_2, n)\pi^* g + \overline{\partial}\partial(\log |z|^2 \cdot \pi^*\rho)$ is a Kähler metric on $\tilde{M}_p$.*
  (2) *On the annulus $\tilde{U} - \tilde{U}'$, we have $\tilde{g} \geq \pi^* g \geq c_2 g_{\varphi'}$ and $\tilde{g} \leq 2A\pi^* g$.*
  (3) *There exists $C(c_1, c_2, n) > 0$ such that with respect to the pull–back coordinate on $\tilde{U} - \tilde{U}'$, we have $\|\tilde{g}\|_{C^2} \leq C(c_1, c_2, n)$.*

**Proof**: We analyze the metric $\tilde{g}$ on the three sets $\tilde{U}'$, $\tilde{M}_p - \tilde{U}$, and $\tilde{U} - \tilde{U}'$ separately.

(1) On $\tilde{U}'$, $\pi^*\rho = 1$, so $\tilde{g} = A\pi^* g + \overline{\partial}\partial \log |z|^2$. Since $\overline{\partial}\partial \log |z|^2 = \overline{\partial}\partial \log |l|^2$ is the Fubini-Study metric on $\mathbb{P}^{n-1}$, the metric $\tilde{g} = A\pi^* g + \overline{\partial}\partial \log |l|^2$, considered as a product metric on $U' \times \mathbb{P}^{n-1}$, is positive definite. Since the metric $\tilde{g}$ on $\tilde{U}'$ can be thought of as the induced metric from the product metric, it is positive definite for any $A > 0$.

(2) On $\tilde{M}_p - \tilde{U}$, $\pi^*\rho = 0$, so $\tilde{g} = A\pi^* g$. Since $\pi$ is biholomorphic on $\tilde{M}_p - \tilde{U}$ and $g$ is positive definite on $M$, $\tilde{g}$ is positive definite on $\tilde{M}_p - \tilde{U}$ for any $A > 0$.

(3) Only on the annulus $\tilde{U} - \tilde{U}'$ do we need a restriction on $A$. Since $\pi$ is biholomorphic on $\tilde{U} - \tilde{U}'$, we have

$$\begin{aligned}
\tilde{g} &= A\pi^* g + \overline{\partial}\partial(\log |z|^2) \cdot \pi^*\rho + \log |z|^2 \cdot \pi^*\overline{\partial}\partial\rho \\
&\quad + \overline{\partial}\log |z|^2 \wedge \pi^*\partial\rho + \partial \log |z|^2 \wedge \pi^*\overline{\partial}\rho \\
&\overset{def}{=} A\pi^* g + h.
\end{aligned}$$

Denote $\frac{\partial}{\partial z'_i}$ by $\partial_{i'}$, where $(\tilde{U} - \tilde{U}', \{z'_1, \cdots, z'_n\})$ is the pull–back coordinates. We have

$$|\partial_{i'} \log |z|^2| = \left|\frac{\overline{z_i}}{|z|^2}\right| \leq 4, \qquad |\partial_{\overline{i'}} \log |z|^2| = \left|\frac{z_i}{|z|^2}\right| \leq 4,$$

$$|\partial_{i'} \partial_{\overline{j'}} \log |z|^2| = \left|\frac{\delta_{ij}|z|^2 - \overline{z_i} z_j}{|z|^4}\right| \leq 32,$$

$$|\pi^* \partial_i \rho| \leq C(n), \quad |\pi^* \partial_{\overline{i}} \rho| \leq C(n), |\pi^* \partial_i \partial_{\overline{j}} \rho| \leq C(n), \quad \text{etc.}$$

Thus there exists a constant $C(n) > 0$ such that

$$|h_{i\overline{j}}|_{C^2} \stackrel{def}{=} |h(\partial_{i'}, \partial_{j'})|_{C^2} \leq C(n).$$

For any vector $\tilde{v} = \sum_i a^i \partial_{i'} \in T_x^{1,0} \tilde{M}_p$, where $x \in \tilde{U} - \tilde{U}'$,

$$|h(\tilde{v}, \tilde{v})| = |\sum_{i,j} a^i \overline{a^j} h_{i\overline{j}}| \leq C(n) \sum_i |a^i|^2 = C(n) g_{\varphi'}(\tilde{v}, \tilde{v}).$$

Since $\varphi' = \varphi \cdot \pi$, we have $\pi_* \partial_{i'} = \partial_i$ for $i = 1, \ldots, n$, and $\pi_* \tilde{v} = \sum_i a^i \partial_i$. So on $\tilde{U} - \tilde{U}'$, we have $g_{\varphi'}(\tilde{v}, \tilde{v}) = g_{\varphi}(\pi_* \tilde{v}, \pi_* \tilde{v})$. Hence

$$\begin{aligned}\tilde{g}(\tilde{v}, \tilde{v}) &\geq (A-1) g(\pi_* \tilde{v}, \pi_* \tilde{v}) - C(n) g_{\varphi'}(\tilde{v}, \tilde{v}) + g(\pi_* \tilde{v}, \pi_* \tilde{v}) \\ &\geq ((A-1)c_2 - C(n)) g_{\varphi'}(\tilde{v}, \tilde{v}) + g(\pi_* \tilde{v}, \pi_* \tilde{v})\end{aligned}$$

For $A = A(c_2, n) \geq \frac{C(n)}{c_2} + 1$, $\tilde{g}$ is positive definite and hence Kähler, which proves the first statement in Lemma 34. Moreover, we get $\tilde{g}(\tilde{v}, \tilde{v}) \geq g(\pi_* \tilde{v}, \pi_* \tilde{v}) \geq c_2 g_{\varphi'}(\tilde{v}, \tilde{v})$, and

$$\begin{aligned}\tilde{g}(\tilde{v}, \tilde{v}) &\leq A \pi^* g(\tilde{v}, \tilde{v}) + C(n) g_{\varphi}(\pi_* \tilde{v}, \pi_* \tilde{v}) \\ &\leq A \pi^* g(\tilde{v}, \tilde{v}) + \frac{C(n)}{c_2} g(\pi_* \tilde{v}, \pi_* \tilde{v}) \\ &\leq 2A g(\pi_* \tilde{v}, \pi_* \tilde{v}),\end{aligned}$$

which proves the second statement. For the third statement, with respect to the pull–back coordinates on $\tilde{U} - \tilde{U}'$, we have

$$\begin{aligned}\tilde{g}_{i\overline{j}} &= A g_{i\overline{j}} + h_{i\overline{j}} \\ |\tilde{g}_{i\overline{j}}|_{C^2} &\leq A(c_2, n) c_1 + C(n) \leq C(c_1, c_2, n).\end{aligned}$$

QED

From now on, let $A$ be chosen as in Lemma 34 so that we have a fixed Kähler metric $\tilde{g}$ on $\tilde{M}_p$. We will compute an upper bound for the diameter of $\tilde{M}_p$ with respect to $\tilde{g}$.

**Proposition 35.** *Let $M$ be a compact Hermitian manifold with $Diam(M) \leq D_0$. Let $\tilde{M}_p$, $(U, \varphi)$ and $A$ be as in the Lemma 34. There exists a constant $C(n, c_1, c_2, D_0) > 0$ such that*

$$Diam_{\tilde{g}}(\tilde{M}_p) \leq C(n, c_1, c_2, D_0),$$

**Proof**: An easy triangle inequality argument shows that

$$\mathrm{Diam}(\tilde{M}_p) \leq 2A\mathrm{Diam}(M) + \mathrm{Diam}(\tilde{U}),$$

so we need to estimate $\mathrm{Diam}(\tilde{U})$. Let $\tilde{x}, \tilde{y}$ be any two points in $\tilde{U}$, and let $x = \pi(\tilde{x})$, $y = \pi(\tilde{y})$. We have three cases:

(1) Assume $x = y = p$. Then $\tilde{x}$ and $\tilde{y}$ lie on the divisor $S = \{0\} \times \mathbb{P}^{n-1}$. Since $\tilde{g}|_S$ is the usual Fubini-Study metric on $\mathbb{P}^{n-1}$, $\mathrm{Dist}(\tilde{x}, \tilde{y}) \leq \mathrm{Diam}_{FS}(\mathbb{P}^{n-1})$, where $\mathrm{Diam}_{FS}(\mathbb{P}^{n-1})$ is the diameter of $\mathbb{P}^{n-1}$ with respect to the Fubini-Study metric.

(2) Assume $x = p$, $y \neq p$. Let

$$\begin{aligned} \tilde{y} &= (z_1, \ldots, z_n, [l_1, \ldots, l_n]) \in \tilde{U} \subset U \times \mathbb{P}^{n-1} \\ \tilde{y}_p &= (0, \ldots, 0, [l_1, \ldots, l_n]) \in \{0\} \times \mathbb{P}^{n-1} \end{aligned}$$

Let $r(t) = (tz_1, \ldots, tz_n, [l_1, \ldots, l_n]) \in U \times \mathbb{P}^{n-1}$, for $0 \leq t \leq 1$. It is easy to see that $r(t) \in \tilde{U}$ and that $r(t)$ is a smooth curve joining $\tilde{y}$ to $\tilde{y}_p$. We have

$$r'(t) = z_1 \partial_1 + \cdots + z_n \partial_n = z_1 \partial_{1'} + \cdots + z_n \partial_{n'}.$$

We will estimate $|r'(t)|$ on $\tilde{U}'$ and $\tilde{U} - \tilde{U}'$ separately:

(A) On $\tilde{U}'$, since $\tilde{g}$ is the induced metric from the product metric, we have

$$|r'(t)|^2 = Ag(\sum_i z_i \partial_i, \sum_j z_j \partial_j) = A \sum_{i,j} z_i \overline{z_j} g_{i\bar{j}} \leq An^2 c_1.$$

(B) On the annulus $\tilde{U} - \tilde{U}'$, we have

$$|r'(t)|^2 = \tilde{g}(\sum_i z_i \partial_{i'}, \sum_j z_j \partial_{j'}) = \sum_{ij} z_i \overline{z_j} \tilde{g}_{i\bar{j}} \leq C(c_1, c_2, n),$$

where the last inequality follows from the third statement in Lemma 34. Hence by (A) and (B) there exists a constant $C(c_1, c_2, n) > 0$ such that $|r'(t)| \leq C(c_1, c_2, n)$. Thus $\mathrm{Dist}(\tilde{y}, \tilde{y}_p) \leq \int_0^1 |r'(t)| dt \leq C(c_1, c_2, n)$, and

$$\mathrm{Dist}(\tilde{x}, \tilde{y}) \leq \mathrm{Dist}(\tilde{y}, \tilde{y}_p) + \mathrm{Dist}(\tilde{x}, \tilde{y}_p) \leq C(c_1, c_2, n) + \mathrm{Diam}_{FS}(\mathbb{P}^{n-1}).$$

(3) Assume $x \neq p$, $y \neq p$. Define $\tilde{x}_p$ analogously to $\tilde{y}_p$. The triangle inequality reduces this case to Cases (1) and (2). Combining (1),(2) and (3), we get the desired result. QED
Before we estimate the curvature of $\tilde{M}_p$, we prove a preparatory lemma.

**Lemma 36.** *On any Riemannian n-manifold $N$, if $(U, \varphi)$ is a local coordinate such that there exist constants $\mu_1, \mu_2$ satisfying (1) $g \geq \mu_1 g_\varphi$, (2) $\|g_{ij}\|_{C^2} \leq \mu_2$, then the absolute value of the sectional curvature on $U$ is bounded by a constant $C(\mu_1, \mu_2, n) > 0$.*

**Proof**: Let $P = \mathrm{span}\{v, w\}$ be any plane in $T_x N$ for $x \in U$, where $v = \sum_i a^i \partial_i$ and $w = \sum_i b^i \partial_i$ are orthornormal vectors. Since

$$1 = g(v, v) \geq \mu_1 \sum_i |a^i|^2, \quad 1 = g(w, w) \geq \mu_1 \sum_i |b^i|^2,$$

we have

$$|a^i|^2 \leq \frac{1}{\mu_1}, \quad |b^i|^2 \leq \frac{1}{\mu_1}, \text{ for } i = 1, \ldots, n.$$

Hence
$$|K(P)| = |\langle R(v,w)w, v\rangle| = |\sum_{i,j,k,l} a^i b^j a^k b^l R_{ijkl}| \leq \frac{C(n)}{\mu_1} \max_{i,j,k,l}\{|R_{ijkl}|\}.$$

Thus it suffices to obtain an upper bound for $|R_{ijkl}|$ in terms of $\mu_1, \mu_2$ and $n$. Since $R_{ijkl}$ can be expressed as a linear combination of products of $g_{ij}, \partial_k g_{ij}, \partial_k \partial_l g_{ij}$ and $g^{ij}, \partial_k g^{ij}$, it is enough to bound $\|g^{ij}\|_{C^1}$ in terms of $\mu_1, \mu_2$ and $n$. By Cramer's rule, it suffices to give a lower bound of $\text{Det}(g)$ in terms of $\mu_1, \mu_2$ and $n$. Since $g \geq \mu_1 g_\varphi$, i.e. $g(v,v) \geq \mu_1 |v|_{E,\varphi}$ for all $v \in T_x N$, the minimum eigenvalue of the matrix $(g_{ij})$ is at least $\mu_1$. Hence $\text{Det}(g) \geq \mu_1^n$.
QED

**Proposition 37.** *Let $M$ be a compact Hermitian manifold with $\text{Ric}(M) \geq k_1$. Let $\tilde{M}_p$, $(U, \varphi)$ and $A$ be as in the Lemma 34. There exists a constant $C(n, c_1, c_2, k_1)$ such that*
$$\text{Ric}_{\tilde{g}}(\tilde{M}_p) \geq C(n, c_1, c_2, k_1).$$

**Proof:** As before we analyze the curvature on the three sets $\tilde{M}_p - \tilde{U}$, $\tilde{U} - \tilde{U}'$ and $\tilde{U}'$ separately.

(1) On $\tilde{M}_p - \tilde{U}$, since the metric $\tilde{g}$ is $A$ times the pull–back metric of the metric $g$, the Ricci curvature on $\tilde{M}_p - \tilde{U}$ is the same as the Ricci curvature on $M$, and hence is at least $k_1$.

(2) On $\tilde{U} - \tilde{U}'$, we take the pull–back coordinates. From Lemma 36 and (2), (3) in Lemma 34, there exists a constant $C(n, c_1, c_2)$ such that $\text{Ric}_{\tilde{g}}(\tilde{U} - \tilde{U}') \geq C(n, c_1, c_2)$.

(3) On $\tilde{U}'$, we take the new coordinates on the open cover $\{\tilde{U}'_1, \dots, \tilde{U}'_n\}$. By relabeling, it suffices to calculate the Ricci curvature on $\tilde{U}'_1$. On $\mathbb{P}^{n-1}$, we have an open cover $\{V_1, \dots, V_n\}$ where

(38) $$V_i = \{[l_1, \dots, l_n] \mid l_i \neq 0 \text{ and } |\frac{l_j}{l_i}| < 10 \text{ if } j \neq i\} \text{ for } i = 1, \dots, n.$$

On $V_1$, we have the standard coordinates $(\tilde{l}) = (\tilde{l}_2, \tilde{l}_3, \dots, \tilde{l}_n)$ where $\tilde{l}_j = \frac{l_j}{l_1}$ for $j = 2, 3, \dots, n$. Since a point in $\tilde{U}'_1$ has two sets of coordinates related by
$$(\tilde{z}_1, \cdots, \tilde{z}_n) \in \tilde{U}'_1 \longleftrightarrow (\tilde{z}_1, \tilde{z}_1 \tilde{z}_2, \cdots, \tilde{z}_1 \tilde{z}_n, [1, \tilde{z}_2, \cdots, \tilde{z}_n]) \in U \times \mathbb{P}^{n-1},$$
we can express $\partial_{\tilde{k}} \stackrel{def}{=} \frac{\partial}{\partial \tilde{z}_k}$ as a tangent vector $U \times \mathbb{P}^{n-1}$

(39) $$\partial_{\tilde{1}} = \partial_1 + \sum_{j=2}^n \frac{l_j}{l_1} \partial_j; \qquad \partial_{\tilde{j}} = z_1 \partial_j + \frac{\partial}{\partial \tilde{l}_j}, \text{ for } j = 2, 3, \dots, n.$$

Thus $\tilde{g}_{kl} = \tilde{g}(\partial_{\tilde{k}}, \partial_{\tilde{l}})$ is given by
$$\tilde{g}_{jk} = Az_1^2 g_{jk} + g_{jk}^{FS}, \text{ for } j \neq 1, k \neq 1,$$
$$\tilde{g}_{1j} = Az_1 g_{1j} + A\sum_{k=2}^n z_k g_{jk}, \text{ for } j \neq 1,$$
$$\tilde{g}_{11} = Ag_{11} + 2A\sum_{j=2}^n \frac{l_j}{l_1} g_{1j} + \sum_{j,k=2}^n \frac{l_j l_k}{l_1^2} g_{jk},$$

where $g^{FS}_{jk} = g_{FS}(\frac{\partial}{\partial \tilde{l}_j}, \frac{\partial}{\partial \tilde{l}_k})$ denotes the Fubini-Study metric. Since $|z_j| \leq \frac{1}{2}, |\frac{l_j}{l_1}| \leq 10, \|g_{ij}\|_{C^2} \leq c_1$ and $\|g^{FS}_{jk}\|_{C^2}$ is a fixed constant, it is easy to check that there exists a constant $C(c_1, c_2, n) > 0$ such that $\|\tilde{g}_{ij}\|_{C^2} \leq C(c_1, c_2, n)$. In the rest of the proof, we will show that there exists a constant $C(c_2, n) > 0$ such that $\tilde{g} \geq C(c_2, n)g_{\tilde{\varphi}_1}$. By Lemma 36, this will finish the proof. We can write every tangent vector on $\tilde{U}' \subset U \times \mathbb{P}^{n-1}$ in the form $\tilde{v} + \tilde{w}$, where $\tilde{v} \in TU'$ and $\tilde{w} \in T\mathbb{P}^{n-1}$. We have

$$\tilde{g}(\tilde{v} + \tilde{w}, \tilde{v} + \tilde{w}) = Ag(\pi_*\tilde{v}, \pi_*\tilde{v}) + g_{FS}(\tilde{w}, \tilde{w}).$$

Since $\overline{V_1}$ is compact, there exists a fixed constant $C_{FS} > 0$ such that $g_{FS} \geq C_{FS} g_{\{\tilde{l}\}}$. Hence we have

$$\begin{aligned}
\tilde{g}(\tilde{v} + \tilde{w}, \tilde{v} + \tilde{w}) &\geq Ac_2 g_\varphi(\pi_*\tilde{v}, \pi_*\tilde{v}) + C_{FS} g_{\{\tilde{l}\}}(\tilde{w}, \tilde{w}) \\
&\geq \min(Ac_2, C_{FS})(g_\varphi(\pi_*\tilde{v}, \pi_*\tilde{v}) + g_{\{\tilde{l}\}}(\tilde{w}, \tilde{w})).
\end{aligned}$$

Let $\tilde{v} + \tilde{w} = \sum_j a^j \partial_{\tilde{j}}$. Then $g_{\tilde{\varphi}_1}(\tilde{v} + \tilde{w}, \tilde{v} + \tilde{w}) = \sum_j |a^j|^2$. By (39) we have

$$\tilde{v} + \tilde{w} = a^1 \partial_1 + \sum_{j=2}^n (a^j z_1 + a^1 \frac{l_j}{l_1})\partial_j + \sum_{j=2}^n a^j \frac{\partial}{\partial \tilde{l}_j}.$$

Hence

$$g_\varphi(\pi_*\tilde{v}, \pi_*\tilde{v}) = |a^1|^2 + \sum_{j=2}^n |(a^j z_1 + a^1 \frac{l_j}{l_1})|^2 \text{ and } g_{\{\tilde{l}\}}(\tilde{w}, \tilde{w}) = \sum_{j=2}^n |a^j|^2,$$

so $g_\varphi(\pi_*\tilde{v}, \pi_*\tilde{v}) + g_{\{\tilde{l}\}}(\tilde{w}, \tilde{w}) \geq g_{\tilde{\varphi}_1}(\tilde{v} + \tilde{w}, \tilde{v} + \tilde{w})$. Thus $\tilde{g} \geq \min(Ac_2, C_{FS})g_{\tilde{\varphi}_1}$. QED

## 6. Generalization of the Kodaria Embedding Theorem for Kähler Manifolds

We first consider the curvature of the holomorphic line bundle $L_p$ associated to the divisor $S$ on $\tilde{M}_p$. Let $\sigma$ denote the projection $\sigma : \tilde{U} \longrightarrow \mathbb{P}^{n-1}$, $\sigma(z, [l]) = [l]$ and let $L_p|_{\tilde{U}}$ denote the restriction of the line bundle $L_p$ to $\tilde{U}$. Then $L_p^*|_{\tilde{U}} = \sigma^* H$, where $H$ is the hyperplane section bundle over $\mathbb{P}^{n-1}$. On the bundle $H \to \mathbb{P}^{n-1}$, we have a natural metric $h_0$ such that the curvature form $\Omega_{FS}$ is the Kähler form of the Fubini-Study metric on $\mathbb{P}^{n-1}$. We equip $L_p^*|_{\tilde{U}}$ with the metric $h_1 = \sigma^* h_0$. Since $L_p^*|_{\tilde{M}_p - \tilde{U}'}$ is trivial, we can equip it with a constant metric $h_2$ with respect to some fixed trivialization. We will fix $h_2$ during the proof of Proposition 40. With $\rho$ chosen as before, we see that $h = \pi^* \rho h_1 + (1 - \pi^* \rho)h_2$ defines a metric on $L_p^*$.

**Proposition 40.** *Let $\tilde{M}_p$, $(U, \varphi)$, $\tilde{g}$ be as in Lemma 34. Then the curvature of $L_p^*$ with respect to the Hermitian metric $h$ above is bounded by a constant $C(c_1, c_2, n) > 0$.*

**Proof:** As before we have three cases:
(1) On $\tilde{M}_p - \tilde{U}$, $h = h_2$ is a constant metric, so $\Omega(L_p^*|_{\tilde{M}_p - \tilde{U}})$ is identically zero.
(2) On $\tilde{U}'$, since $h = h_1 = \sigma^* h_0$, we have $\Omega(L_p^*|_{\tilde{U}'}) = \sigma^* \Omega_{FS}$. For any tangent vector $\tilde{v} \in T^{1,0}\tilde{U}'$,

$$|\Omega(L_p^*|_{\tilde{U}'})(\tilde{v}, \tilde{v})| = |\Omega_{FS}(\sigma_*\tilde{v}, \sigma_*\tilde{v})| = g_{FS}(\sigma_*\tilde{v}, \sigma_*\tilde{v}) \leq \tilde{g}(\tilde{v}, \tilde{v}).$$

Hence $\Omega(L_p^*|_{\tilde{U}'})$ is bounded by the constant 1.
(3) For $\tilde{U} - \tilde{U}'$, we will find a nonwhere zero holomorphic section $s$ of $L_p^*$ and compute

x

$\Omega(L_p^*) = -\partial\overline{\partial}\log h$. On $\mathbb{P}^{n-1}$ we choose the open cover $\{V_1, \ldots, V_n\}$ as in (38). There exists a nowhere vanishing holomorphic section $f_k$ of the hyperplane section bundle $H$ over $\overline{V_k}$ with the transition functions $\frac{l_k}{l_j}$ over $V_j \cap V_k$. Since $\overline{V_k}$ is compact, there exist positive constants $\tau_1, \tau_2$ such that for $k = 1, \ldots, n$, we have

$$\tau_1 \le h_0(f_k, f_k) \le \tau_2, \quad \left|\frac{\partial}{\partial \tilde{l}_i} h_0(f_k, f_k)\right| \le \tau_2, \quad \text{and} \quad \left|\frac{\partial^2}{\partial \tilde{l}_i \partial \overline{\tilde{l}_j}} h_0(f_k, f_k)\right| \le \tau_2 \text{ on } \overline{V_k},$$

where $\{\tilde{l}_1, \ldots, 1, \ldots, \tilde{l}_n\}$ are the natural coordinates on $V_k$. We have an open cover $\{\tilde{U}_k = \sigma^{-1}V_k\}$ of $\tilde{U}$ as before. Since $L_p^*|_{\tilde{U}} = \sigma^*H$, $e_k \stackrel{def}{=} \sigma^* f_k = f_k \circ \sigma$ is a local frame of $L_p^*$ over $\tilde{U}_k$, and the transition functions of $L_p^*$ are given by $\frac{l_k}{l_j} = \frac{z_k}{z_j}$ over $\tilde{U}_j \cap \tilde{U}_k$. We claim that $|z_k|$ is bounded away from zero on $(\tilde{U} - \tilde{U}') \cap \tilde{U}_k$. In fact, since $z_j = z_k \frac{l_j}{l_k}$, for all $j$, we have

$$|z|^2 = |z_1|^2 + \cdots + |z_k|^2 = |z_k|^2 (|\frac{l_1}{l_k}|^2 + \cdots + |\frac{l_n}{l_k}|^2) \le |z_k|^2 100n.$$

Since $|z| \ge \frac{1}{2}$ on $\tilde{U} - \tilde{U}'$, we have $|z_k|^2 \ge \frac{1}{400n} = C(n)$, which proves the claim.

Since $z_k$ is never zero on $(\tilde{U} - \tilde{U}') \cap \tilde{U}_k$, we have a global section $s$ of $L_p^*$ on $\tilde{U} - \tilde{U}'$ defined by

$$s = \frac{1}{z_k} e_k \quad \text{over} \quad (\tilde{U} - \tilde{U}') \cap \tilde{U}_k.$$

Choose $h_2$ to be a fixed constant metric such that $h_2(s, s) = 1$. Thus on $(\tilde{U} - \tilde{U}') \cap \tilde{U}_k$, we have

(41) $$h(s,s) = \pi^*\rho h_1(s,s) + (1 - \pi^*\rho) = \pi^*\rho \frac{1}{|z_k|^2} h_0(f_k, f_k) \circ \sigma + (1 - \pi^*\rho)$$

It is enough to compute the upper bound of $|\Omega(L_p^*)|$ on $(\tilde{U} - \tilde{U}') \cap \tilde{U}_1$ only. We will use the pull–back coordinate $((\tilde{U} - \tilde{U}') \cap \tilde{U}_1, \varphi_1') = ((\tilde{U} - \tilde{U}') \cap \tilde{U}_1, z_1', \ldots, z_n')$, where $z_j' = z_j$ for $j = 1, \ldots, n$. Since $\tilde{g} \ge c_2 g_{\varphi_1'}$ it is enough to bound $|\Omega(L_p^*)(\partial_{i'}, \partial_{j'})|$ for all $i, j$ in terms of $n, c_1, c_2$. Since a point in $(\tilde{U} - \tilde{U}') \cap \tilde{U}_1$ has two sets of coordinates related by

$$(z_1', \cdots, z_n') \in (\tilde{U} - \tilde{U}') \cap \tilde{U}_1 \longleftrightarrow (z_1', \cdots, z_n', [1, \frac{z_2'}{z_1'}, \cdots, \frac{z_n'}{z_1'}]) \in U \times \mathbb{P}^{n-1},$$

we can compute $\partial_{i'}$ as a vector in $TU \oplus T\mathbb{P}^{n-1}$

$$\partial_{1'} = \partial_1 + \sum_{j=2}^{n} (-\frac{l_j}{l_1} \frac{1}{z_1}) \frac{\partial}{\partial \tilde{l}_j}; \quad \partial_{k'} = \partial_k + \frac{1}{z_1} \frac{\partial}{\partial \tilde{l}_k}, \quad \text{for} \quad k = 2, \ldots, n,$$

where $\tilde{l}_j = \frac{l_j}{l_1}$ for $j = 2, \ldots, n$. Let $h = h(s,s)$ and $h_1 = h_1(s,s)$. Then (41) becomes $h = \pi^*\rho h_1 + (1 - \pi^*\rho)$. An easy computation shows that on $(\tilde{U} - \tilde{U}') \cap \tilde{U}_1$, we have

$$\begin{aligned}\Omega(L_p^*) &= -\partial\overline{\partial}\log h \\ &= -\frac{1}{h^2}(\overline{\partial}\pi^*\rho \, (h_1 - 1) + \pi^*\rho \, \overline{\partial} h_1) \wedge (\partial\pi^*\rho \, (h_1 - 1) + \pi^*\rho \, \partial h_1) \\ &\quad + \frac{1}{h}(\overline{\partial}\partial\pi^*\rho \, (h_1 - 1) + \pi^*\rho \, \overline{\partial}\partial h_1 + \overline{\partial}\pi^*\rho \wedge \partial h_1 + \partial\pi^*\rho \wedge \overline{\partial} h_1).\end{aligned}$$

Since $|h_1| \geq 4\tau_1$, we have $h \geq \min\{4\tau_1, 1\}$ and $\frac{1}{h}$ is bounded from above. From

$$\begin{aligned}\partial_{1'} h_1 &= \partial_{1'}\left(\frac{1}{|z_1|^2} h_0(f_1, f_1) \circ \sigma\right) \\ &= \frac{\overline{z_1}}{|z_1|^4} h_0(f_1, f_1) \circ \sigma + \frac{1}{|z_1|^2} \sum_{j=2}^{n} (-\frac{l_j}{l_1}\frac{1}{z_1}) \frac{\partial}{\partial \tilde{l}_j} h_0(f_1, f_1)\end{aligned}$$

and the bounds $|z_1| \geq C(n)$, $|\frac{l_j}{l_1}| < 10$ for $j = 2, \ldots, n$, $h_0(f_1, f_1) \leq \tau_2$, $|\frac{\partial}{\partial \tilde{l}_j} h_0(f_1, f_1)| \leq \tau_2$, we see $|\partial_{1'} h_1|$ is bounded from above by a constant depending only on $\tau_1, \tau_2, n$. Similarly, we can obtain bounds on $\partial_{\overline{1'}} h_1, \partial_{i'} h_1, \partial_{\overline{i'}} h_1, \partial_{i'} \partial_{\overline{j'}} h_1$, etc. Since $|\partial_i \pi^* \rho|$, $|\partial_{\overline{i}} \pi^* \rho|$, and $|\partial_i \partial_{\overline{j}} \pi^* \rho|$ are bounded from above by a constant depending only on $n$, we get an upper bound for $|\Omega(L_p^*)(\partial_{i'}, \partial_{j'})|$ for all $i, j$. QED

We can now prove a generalization of the Kodaria embedding theorem for Kähler manifolds.

**Theorem 42.** *Let $\mathcal{M} = \mathcal{M}(n, D_0, V_0, k_1, k_2, c_1, c_2)$ be the collection of pairs $(M, F)$ where*

(1) *$M$ is a compact Kähler manifold of complex dimension $n$ with $\mathrm{Vol}(M) \geq V_0$, $\mathrm{Diam}(M) \leq D_0$, $|\mathcal{R}ic(M)| \leq k_1$;*
(2) *On $M$ there exists an open set $U$ with volume less than $\frac{1}{2} V_0$ and a local coordinate $(U, \varphi)$ of class $\mathcal{C}(c_1, c_2)$;*
(3) *$F$ is a Hermitian holomorphic line bundle over $M$ with curvature $|\Omega(F)| \leq k_2$.*

*Given $\alpha > 0$, there exists $c = c(\mathcal{M}, \alpha) > 0$ such that for any $(M, F) \in \mathcal{M}$ if $\Omega(F) > \alpha$ except on a set in $M - U$ of volume less than $c$, then $M$ is projective algebraic.*

**Proof:** Let $p \in U$ have $\varphi(p) = 0$. Let $\tilde{M}_p$ be the blow–up manifold of $M$, constructed as in §5 from the local coordinate $(U, \varphi)$ centered at $p$. Let $\tilde{g}$ be the Kähler metric constructed in Lemma 34. From Propositions 35 and 37, we know that there exist constants

$$\tilde{D}_0 = \tilde{D}_0(n, c_1, c_2, D_0) > 0, \quad \tilde{k}_1 = \tilde{k}_1(n, c_1, c_2, k_1)$$

such that

$$\mathrm{Diam}(\tilde{M}_p) \leq \tilde{D}_0, \quad \mathcal{R}ic(\tilde{M}_p) \geq \tilde{k}_1.$$

We have $\mathrm{Vol}(M - U) > \frac{1}{2} V_0$. Since $\tilde{g} = A\pi^* g$ on $\tilde{M}_p - \tilde{U}$, we have

$$\mathrm{Vol}(\tilde{M}_p) > \mathrm{Vol}(\tilde{M}_p - \tilde{U}) = A^{\frac{n}{2}} \mathrm{Vol}(M - U) > \frac{1}{2} A^{\frac{n}{2}} V_0 \stackrel{def}{=} \tilde{V}_0(n, c_1, c_2, V_0).$$

From Corollary 32 and Lemma 33, we know that it suffices to prove the vanishing of $H^1(\tilde{M}_p, \pi^* E \otimes L_p^{*\mu})$ for $\mu = 1, 2$ and some holomorphic line bundle $E$. In order to apply our generalization of the Kodaira vanishing theorem (Theorem 25) to the bundle $\pi^* E \otimes L_p^{*\mu}$, we need to compute the lowest eigenvalue of $\Omega(\pi^* E \otimes L_p^{*\mu} \otimes K_{\tilde{M}_p}^*)$ and a bound for $\Omega(\pi^* E \otimes L_p^{*\mu})$. Since $K_{\tilde{M}_p}^* = \pi^*(K_M^*) \otimes L_p^{*(n-1)}$, we have $\Omega(\pi^* E \otimes L_p^{*\mu} \otimes K_{\tilde{M}_p}^*) = \Omega(\pi^*(E \otimes K_M^*) \otimes L_p^{*(n-1+\mu)})$. Choosing $E = F^k \otimes K_M$, where $k$ is an integer to be determined below, we need to compute (a) the lowest eigenvalue of $\Omega(\pi^*(F^k) \otimes L_p^{*(n-1+\mu)})$, and (b) a bound for $\Omega(\pi^*(F^k \otimes K_M) \otimes L_p^{*\mu})$. Denote the lowest eigenvalue of $\Omega$ by $\lambda(\Omega)$.

**Case I:** $\tilde{U} - \tilde{U}'$

On the annulus $\tilde{U} - \tilde{U}'$, by assumption we know $\Omega(F) > \alpha$ on $U$. For any tangent vector $\tilde{v}$ to $\tilde{U} - \tilde{U}'$, we have

$$\pi^*\Omega(F)(\tilde{v}, \tilde{v}) \geq \alpha g(\pi_*\tilde{v}, \pi_*\tilde{v}) \geq \frac{\alpha}{2A}\tilde{g}(\tilde{v}, \tilde{v})$$

by (2) in Lemma 34. By Proposition 40, we know $|\lambda(\Omega(L_p^*))|$ is bounded by a constant $C(n, c_1, c_2) > 0$. Since

$$\begin{aligned}\lambda(\Omega(\pi^*(F^k) \otimes L_p^{*(n-1+\mu)})) &\geq k\lambda(\pi^*\Omega(F)) + (n - 1 + \mu)\lambda(\Omega(L_p^*)) \\ &\geq \frac{k\alpha}{2A} - (n - 1 + \mu)C(n, c_1, c_2),\end{aligned}$$

we can choose $k = k(n, c_1, c_2, \alpha) > 1$ such that

$$(a) \quad \lambda(\Omega(\pi^*(F^k) \otimes L_p^{*(n-1+\mu)})) \geq \frac{\alpha}{2A} \text{ on } \tilde{U} - \tilde{U}'.$$

We now fix such a $k$. Still on the annulus $\tilde{U} - \tilde{U}'$, we have

$$\begin{aligned}\left|\Omega(\pi^*(F^k \otimes K_M) \otimes L_p^{*\mu})(\tilde{v}, \tilde{v})\right| &\leq |k\pi^*\Omega(F)(\tilde{v}, \tilde{v}) + \pi^*\Omega(K_M)(\tilde{v}, \tilde{v}) + \mu\Omega(L_p^*)(\tilde{v}, \tilde{v})| \\ &\leq |k\Omega(F)(\pi_*\tilde{v}, \pi_*\tilde{v}) + \Omega(K_M)(\pi_*\tilde{v}, \pi_*\tilde{v}) + \mu\Omega(L_p^*)(\tilde{v}, \tilde{v})| \\ &\leq (kk_2 + k_1)g(\pi_*\tilde{v}, \pi_*\tilde{v}) + \mu C(n, c_1, c_2)\tilde{g}(\tilde{v}, \tilde{v}) \\ &\leq (kk_2 + k_1 + \mu C(n, c_1, c_2))\tilde{g}(\tilde{v}, \tilde{v}),\end{aligned}$$

which says that

$$(b) \quad \left|\Omega(\pi^*(F^k \otimes K_M) \otimes L_p^{*\mu})\right| \leq kk_2 + k_1 + \mu C(n, c_1, c_2) \text{ on } \tilde{U} - \tilde{U}'.$$

**Case II:** $\tilde{M}_p - \tilde{U}$

Let $B$ be the set where $\Omega(F) \leq \alpha$. For any tangent vector $\tilde{v}$ of type (1,0) to $\tilde{M}_p - \tilde{U} - \pi^{-1}(B)$, we have

$$\Omega(\pi^*(F^k) \otimes L_p^{*(n-1+\mu)})(\tilde{v}, \tilde{v}) = k\Omega(F)(\pi_*\tilde{v}, \pi_*\tilde{v}) \geq k\alpha g(\pi_*\tilde{v}, \pi_*\tilde{v}) \geq \frac{k\alpha}{A}\tilde{g}(\tilde{v}, \tilde{v}),$$

which says that

$$(a) \quad \lambda(\Omega(\pi^*(F^k) \otimes L_p^{*(n-1+\mu)})) \geq \frac{\alpha}{A} \text{ on } \tilde{M}_p - \tilde{U} - \pi^{-1}(B).$$

On the other hand, on all of $\tilde{M}_p - \tilde{U}$,

$$\begin{aligned}\left|\Omega(\pi^*(F^k \otimes K_M) \otimes L_p^{*\mu})(\tilde{v}, \tilde{v})\right| &= |k\Omega(F)(\pi_*\tilde{v}, \pi_*\tilde{v}) + \Omega(K_M)(\pi_*\tilde{v}, \pi_*\tilde{v})| \\ &\leq (kk_2 + k_1)g(\pi_*\tilde{v}, \pi_*\tilde{v}) \leq \frac{kk_2 + k_1}{A}\tilde{g}(\tilde{v}, \tilde{v}),\end{aligned}$$

which says that

$$(b) \quad \left|\Omega(\pi^*(F^k \otimes K_M) \otimes L_p^{*\mu})\right| \leq \frac{kk_2 + k_1}{A} \text{ on } \tilde{M}_p - \tilde{U}.$$

**Case III:** $\tilde{U}'$

On $\tilde{U}'$, we have $\Omega(L_p^*) = \sigma^*\Omega_{FS}$. Any tangent vector to $\tilde{U}'$ can be written in form $\tilde{v} + \tilde{w}$, where $\tilde{v} \in TU'$ and $\tilde{w} \in T\mathbb{P}^{n-1}$. We have

$$\begin{aligned}
\Omega(\pi^*(F^k) \otimes L_p^{*(n-1+\mu)})(\tilde{v}+\tilde{w}, \tilde{v}+\tilde{w}) &= k\Omega(F)(\pi_*\tilde{v}, \pi_*\tilde{v}) + (n-1+\mu)\Omega_{FS}(\sigma_*\tilde{w}, \sigma_*\tilde{w}) \\
&\geq k\alpha g(\pi_*\tilde{v}, \pi_*\tilde{v}) + (n-1+\mu)g_{FS}(\sigma_*\tilde{w}, \sigma_*\tilde{w}) \\
&\geq \min(\frac{k\alpha}{A}, n-1+\mu)(Ag(\pi_*\tilde{v}, \pi_*\tilde{v}) + g_{FS}(\sigma_*\tilde{w}, \sigma_*\tilde{w})) \\
&\geq \min(\frac{k\alpha}{A}, n-1+\mu)\tilde{g}(\tilde{v}+\tilde{w}, \tilde{v}+\tilde{w}),
\end{aligned}$$

which says that

$$(a) \quad \lambda(\Omega(\pi^*(F^k) \otimes L_p^{*(n-1+\mu)})) \geq \min(\frac{k\alpha}{A}, n-1+\mu) \text{ on } \tilde{U}'.$$

On the other hand,

$$\begin{aligned}
\left|\Omega(\pi^*(F^k \otimes K_M) \otimes L_p^{*\mu})(\tilde{v}+\tilde{w}, \tilde{v}+\tilde{w})\right| &\leq |k\pi^*\Omega(F)(\tilde{v}, \tilde{v}) + \pi^*\Omega(K_M)(\tilde{v}, \tilde{v})| + |\mu\sigma^*\Omega_{FS}(\tilde{w}, \tilde{w})| \\
&\leq (kk_2 + k_1)g(\pi_*\tilde{v}, \pi_*\tilde{v}) + \mu g_{FS}(\sigma_*\tilde{w}, \sigma_*\tilde{w}) \\
&\leq \max\left(\frac{kk_2 + k_1}{A}, n-1+\mu\right)\tilde{g}(\tilde{v}+\tilde{w}, \tilde{v}+\tilde{w}),
\end{aligned}$$

which says that

$$(b) \quad \left|\Omega(\pi^*(F^k \otimes K_M) \otimes L_p^{*\mu})\right| \leq \max\left(\frac{kk_2 + k_1}{A}, \mu\right) \text{ on } \tilde{U}'.$$

Take $\tilde{\alpha} = \min(\alpha/2A, n-1+\mu) > 0$ and $\tilde{k}_2 = kk_2 + k_1 + \mu C(c_1, c_2, n) + \mu$. So far we have shown that

$$(a) \quad \lambda(\Omega(\pi^*(F^k) \otimes L_p^{*(n-1+\mu)})) > \tilde{\alpha}$$

except on a set $\pi^{-1}(B)$, and (since $A > 1$)

$$(b) \quad \left|\Omega(\pi^*(F^k \otimes K_M) \otimes L_p^{*\mu})\right| \leq \tilde{k}_2$$

on $\tilde{M}_p$. By our generalization of the Kodaria vanishing theorem (Theorem 25), there exists a constant $c(n, \tilde{D}_0, \tilde{V}_0, \tilde{k}_1, \tilde{k}_2, \tilde{\alpha}) > 0$ such that if

$$\text{Vol}(\pi^{-1}(B)) < c(n, \tilde{D}_0, \tilde{V}_0, \tilde{k}_1, \tilde{k}_2, \tilde{\alpha})$$

then $H^1(\tilde{M}_p, \pi^*(F^k \otimes K_M) \otimes L_p^{*\mu})) = 0$, for $\mu = 1, 2$, which implies that $M$ is algebraic projective by Lemma 33 and Corollary 32. Since $\tilde{g} = A\pi^*g$ on $\pi^{-1}(B)$, we have $\text{Vol}(\pi^{-1}(B)) = A^{\frac{n}{2}}\text{Vol}(B)$. Let

$$c(\mathcal{M}, \alpha) = c(n, D_0, V_0, k_1, k_2, c_1, c_2, \alpha) \stackrel{def}{=} \frac{c(n, \tilde{D}_0, \tilde{V}_0, \tilde{k}_1, \tilde{k}_2, \tilde{\alpha})}{A^{\frac{n}{2}}}$$

Therefore, if $\text{Vol}(B) < c(\mathcal{M}, \alpha)$, $M$ is projective algebraic. QED

DEPARTMENT OF MATHEMATICS, BOSTON UNIVERSITY, BOSTON, MASSACHUSETTS 02215
*E-mail address*: ying@math.bu.edu